\documentclass[]{article}

\usepackage[english]{babel}
\usepackage{amsmath}
\usepackage{amsthm}
\usepackage{amsfonts}
\usepackage{epsf}
\usepackage{graphicx}

\numberwithin{theorem}{section}
\numberwithin{definition}{section}

\begin{document}
	
\title{Taxicab Calculus: Trig Derivatives}
\author{Kevin P. Thompson}
	
\date{}
	
\thispagestyle{empty}
\renewcommand\thispagestyle[1]{} 
	
\maketitle
	
\begin{abstract}
The set of trigonometric functions in taxicab geometry is completed and derivatives of all of such functions are explored.
\end{abstract}

\maketitle

\section{Introduction}

In \cite{ThompsonDray, Akca} the foundation was laid for native angles and trigonometry in taxicab geometry. But, each of these initial works was incomplete. In particular, the secant and cosecant functions were not explored. As we quickly review the existing work and tie up these loose ends, some very intriguing characteristics will be apparent in some of the taxicab trigonometric functions. This will lead us into exploring the differentiability of these functions and the expressions of such derivatives. We will begin by briefly reviewing the established functions and developing those that have been missed to this time: secant and cosecant.

\section{Taxicab Trigonometry}

\begin{figure}[t]
\centering
\includegraphics[width=0.6\linewidth]{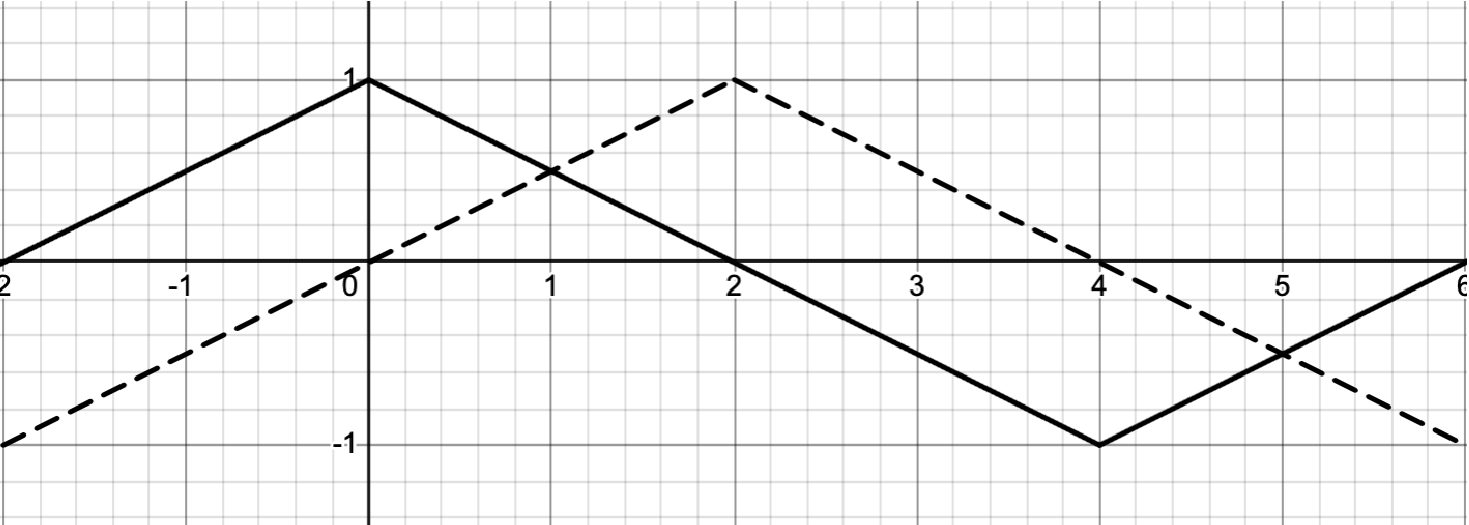}
\caption{Graphs of the taxicab sine (dashed) and cosine (solid) functions.}
\label{sincos}
\end{figure}

Traditionally, taxicab geometry simply used Euclidean angles (see \cite{Krause}), a variant that should be dubbed "modified taxicab geometry." This research will use pure taxicab geometry where angles are natively defined using taxicab distance along the unit taxicab circle. Such angles are measured in t-radians (see \cite{ThompsonDray}). Since the unit taxicab circle has a circumference of 8, this leads to cosine and sine functions with period $8 = 2\pi_t$ where $\pi_t = 4$ is the value for $\pi$ in taxicab geometry (see \cite{Euler}). The definitions of cosine and sine from \cite{ThompsonDray} are given below with their graphs shown in Figure \ref{sincos}.
\[
\cos \theta =\left\{
\begin{array}{l}
1-\frac{1}{2}\theta \;,\;\;\;0\leq \theta <\pi_t \\
\\
-3+\frac{1}{2}\theta \;,\;\;\;\pi_t\leq \theta <2\pi_t
\end{array}
\right.  \;,\;\;\;\;\sin \theta =\left\{
\begin{array}{l}
\frac{1}{2}\theta \;,\;\;\;-\frac{\pi_t}{2}\leq \theta <\frac{\pi_t}{2} \\
\\
2-\frac{1}{2}\theta \;,\;\;\;\frac{\pi_t}{2}\leq \theta <\frac{3\pi_t}{2} \\
\\
-4+\frac{1}{2}\theta \;,\;\;\;\frac{3\pi_t}{2}\leq \theta <\frac{5\pi_t}{2}
\end{array}
\right.
\]

With a bit of effort, closed-form expressions can be invented for these piecewise functions. On the interval $[(k-1)\frac{\pi_t}{2},k\frac{\pi_t}{2})$, the closed-form expressions of sine and cosine are
\[
\cos \theta = \frac{(1+(-1)^k)i^{k-2}(k-1-\frac{1}{2}\theta)+(1+(-1)^{k+1})i^{k-1}(k-\frac{1}{2}\theta) }{2}
\]
\[
\sin \theta = \frac{(1+(-1)^k)i^{k-2}(k-\frac{1}{2}\theta)+(1+(-1)^{k+1})i^{k-1}(1-k+\frac{1}{2}\theta)}{2}
\]

These closed-form expressions are a bit unweildy to deal with. And, if we look forward to where we are headed, one of the sine and cosine functions switches branches in the piecewise definition every $\frac{\pi_t}{2}$ t-radians. So, the expressions we choose to use will be most useful when they are based on $\frac{\pi_t}{2}$ t-radian intervals. This will lend itself nicely when we consider functions that are combinations of both sine and cosine such as the tangent function. On the interval $[(k-1)\frac{\pi_t}{2},k\frac{\pi_t}{2})$, the pseudo closed-form expressions of sine and cosine are
\[
   \cos \theta = \left\{
     \begin{array}{ll}
       i^{k-2}(k-1-\frac{1}{2}\theta), & \text{if }k\text{ is even} \\
       i^{k-1}(k-\frac{1}{2}\theta), & \text{if }k\text{ is odd}
     \end{array}
   \right.
\]
\[
   \sin \theta = \left\{
     \begin{array}{ll}
       i^{k-2}(k-\frac{1}{2}\theta), & \text{if }k\text{ is even} \\
       i^{k-1}(1-k+\frac{1}{2}\theta), & \text{if }k\text{ is odd}
     \end{array}
   \right.
\]
where the imaginary number $i$ is used to create positive and negative terms as needed. As an example, choosing $k=1$ for these gives the expressions $\cos \theta = 1 - \frac{\theta}{2}$ and $\sin \theta = \frac{\theta}{2}$ on the interval $[0,\frac{\pi_t}{2})$ in agreement with our earlier piecewise definition.

If we define the taxicab tangent function as in Euclidean geometry as the ratio of sine and cosine, we again have a piecewise-defined function. Since we will be using this in later analysis, the pseudo closed-form expression of the taxicab tangent function is given below. The graph of this function, shown in Figure \ref{tangent}, has a very similar feel to the Euclidean tangent function with asymptotes at $k\pi_t-\frac{\pi_t}{2}$ t-radians and a period that is half that of sine and cosine. On the interval $[(k-1)\frac{\pi_t}{2},k\frac{\pi_t}{2})$, 
\[
   \tan \theta = \left\{
     \begin{array}{ll}
       \frac{k-\frac{1}{2}\theta }{k-1-\frac{1}{2}\theta }, & \text{if }k\text{ is even} \\
       \frac{1-k+\frac{1}{2}\theta }{k-\frac{1}{2}\theta }, & \text{if }k\text{ is odd}
     \end{array}
   \right.
\]
\begin{figure}[t]
\centering
\includegraphics[width=0.5\linewidth]{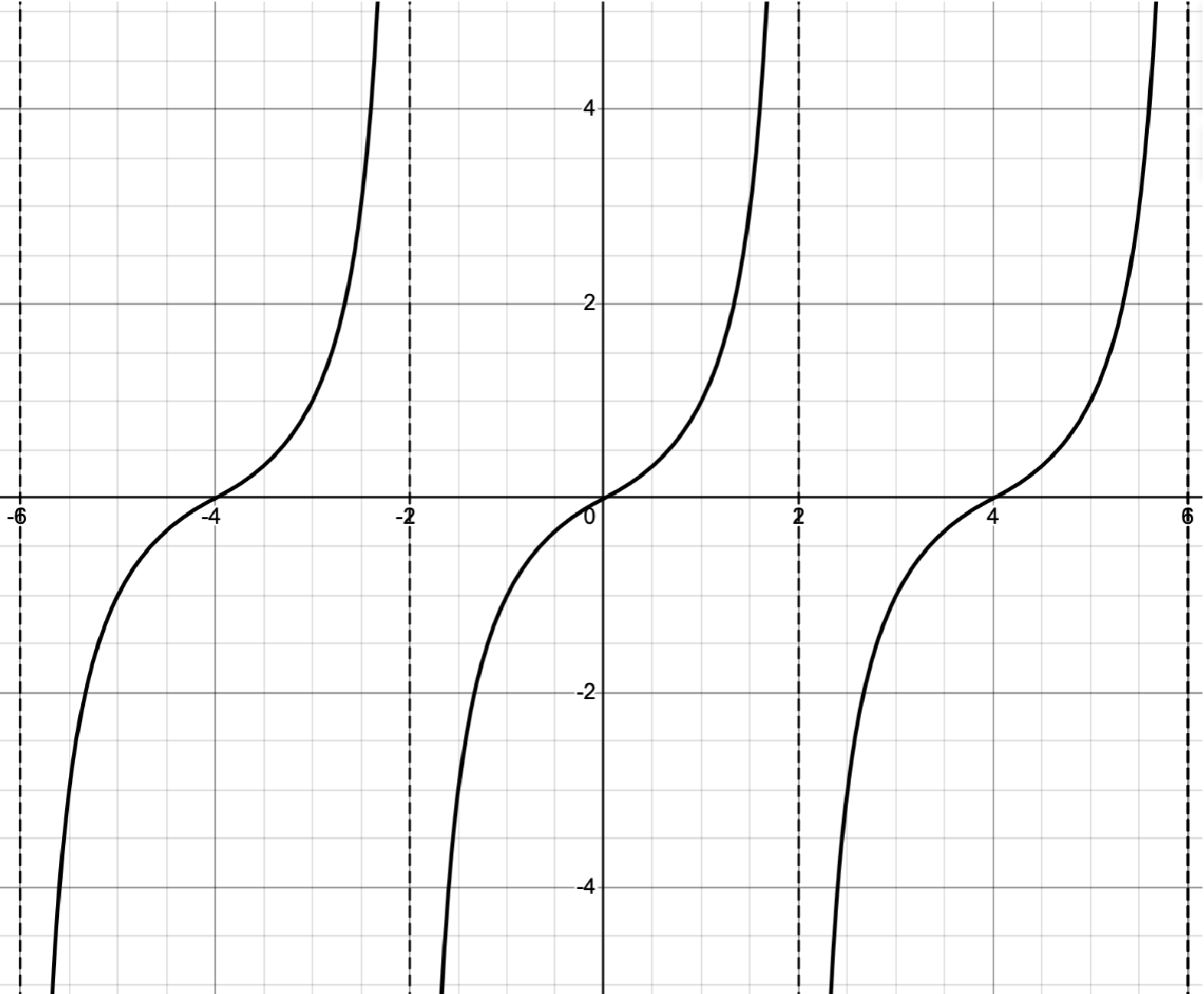}
\caption{Graph of the taxicab tangent function.}
\label{tangent}
\end{figure}

On the surface, the taxicab tangent function appears to be differentiable everywhere between the asymptotes except where cosine is zero. A formal expression for the derivative would confirm this observation.

\begin{figure}[b]
\centering
\includegraphics[width=0.45\linewidth]{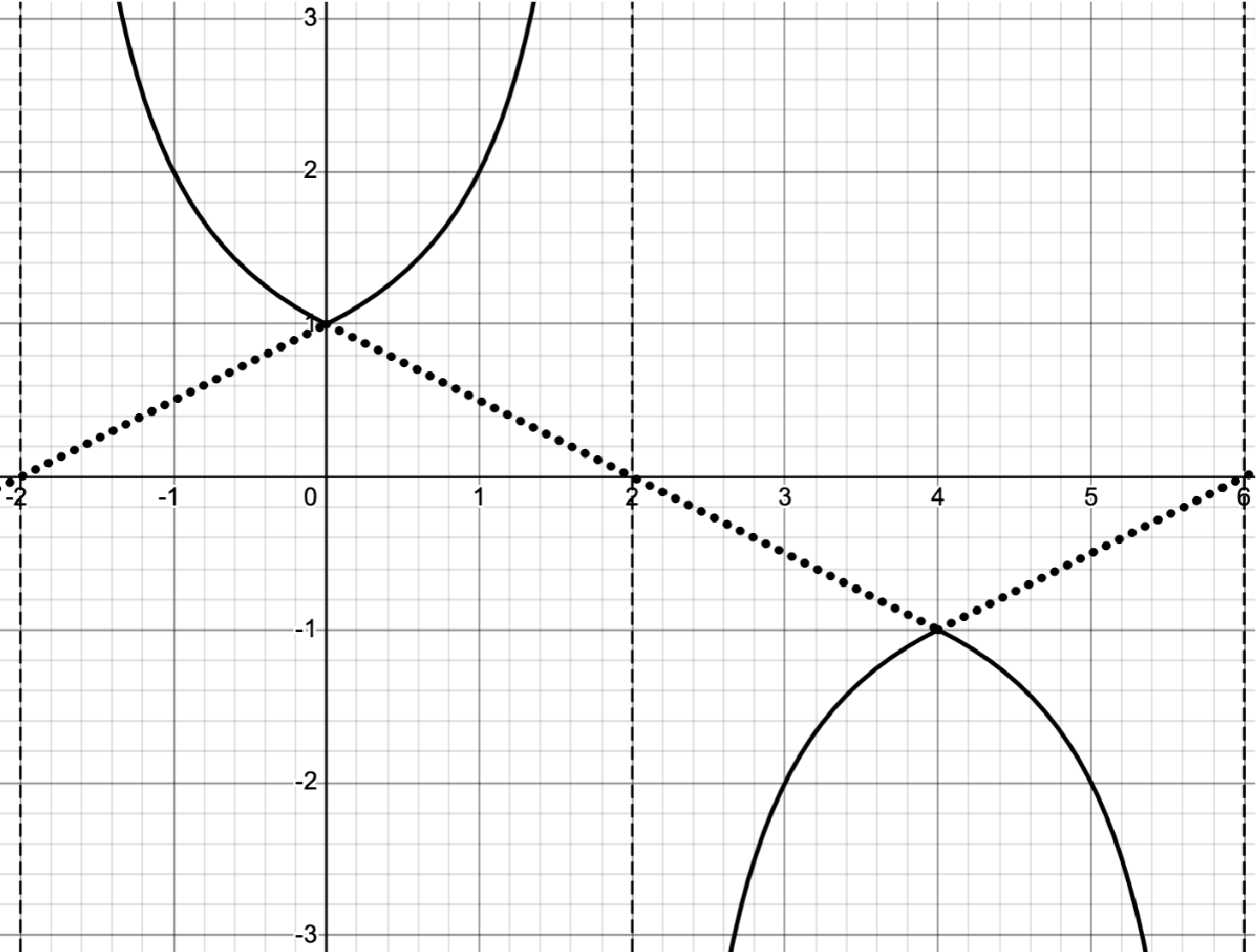}
\caption{Graphs of the taxicab secant (solid) and cosine (dotted) functions.}
\label{secant}
\end{figure}

Since cotangent is simply the reciprocal of tangent, we will leave this to the reader and move to the secant function. Using the reciprocals of the various piecewise branches of the cosine function, we can construct a graph of secant (shown with cosine in Figure \ref{secant}).

The taxicab secant graph exhibits the same characteristics as the Euclidean secant graph: a period equal to that of cosine; asymptotes where cosine is zero; and, a minimum or maximum where the cosine has a maximum or minimum, respectively. As with tangent, the secant graph is curved. But, unlike the Euclidean secant function, the taxicab secant function is not differentiable at its extrema. This is to be expected since this is where cosine is not differentiable.

As with cotangent, we will leave the details of the cosecant function to the reader.

\section{Derivatives of Trig Functions}

Before even beginning this discussion, it is important to consider whether the normal application of calculus is valid in taxicab geometry. A derivative represents the slope of a tangent line. And, slope is a ratio of vertical and horizontal distances. Since the Euclidean and taxicab metrics agree for line segments parallel to the $x$-axis or $y$-axis (see \cite{Thompson}), the concept of slope, and therefore the concept of a derivative, transfers seamlessly to taxicab geometry. (Although we will not consider it here, the concepts of integration and area also transfer seamlessly to taxicab geometry.)

From the definition of the taxicab cosine and sine functions, it is clear that the derivatives of these functions alternate from $\frac{1}{2}$ to $-\frac{1}{2}$ about the extrema. Not only are their derivatives now piecewise constant, the derivative relationship between sine and cosine seen in Euclidean geometry has been lost. Finally, sine and cosine are not everywhere differentiable because of the sharp corners at the extrema. 

If we simply take the definition of tangent as the ratio of sine to cosine, use the closed-form expressions, and differentiate with respect to $\theta$, we obtain the derivative of the taxicab tangent function. As explained above, we will take the proof in two cases: on intervals $[(k-1)\frac{\pi_t}{2},k\frac{\pi_t}{2})$ with $k$ even and odd.
\begin{eqnarray*}
\frac{d}{d\theta}(\tan \theta) & = & \frac{d}{d\theta}\left[ \frac{i^{k-2}(k-\frac{1}{2}\theta)}{i^{k-2}(k-1-\frac{1}{2}\theta)} \right] \\
        & = & \frac{i^{k-2}(k-1-\frac{1}{2}\theta)(-\frac{1}{2}i^{k-2})-i^{k-2}(k-\frac{1}{2}\theta)(-\frac{1}{2}i^{k-2})}{(i^{k-2}(k-1-\frac{1}{2}\theta))^2} \\
        & = & \frac{i^{2k-4}(-\frac{1}{2}k+\frac{1}{2}+\frac{1}{4}\theta +\frac{1}{2}k-\frac{1}{4}\theta)}{(i^{k-2}(k-1-\frac{1}{2}\theta))^2} \\
        & = & \frac{\frac{1}{2}(i^2)^k i^{-4}}{(i^{k-2}(k-1-\frac{1}{2}\theta))^2} \text{, (} k\text{ even)}\\
        & = & \frac{\frac{1}{2}}{(i^{k-2}(k-1-\frac{1}{2}\theta))^2} = \frac{\frac{1}{2}}{\cos^2\theta} = \frac{1}{2}\sec^2\theta
\end{eqnarray*}

\noindent For $k$ odd, picking up at the 4th line, we have
\begin{eqnarray*}
	\frac{d}{d\theta}(\tan \theta) & = & \dots \\
	& = & \frac{\frac{1}{2}(i^2)^k i^{-2}}{(i^{k-1}(k-\frac{1}{2}\theta))^2} \text{, (} k\text{ odd)}\\
	& = & \frac{\frac{1}{2}}{(i^{k-1}(k-\frac{1}{2}\theta))^2} = \frac{\frac{1}{2}}{\cos^2\theta} = \frac{1}{2}\sec^2\theta
\end{eqnarray*}

So, the derivative of tangent is $\frac{1}{2}\sec^2 \theta$, and we can confirm that the taxicab tangent function is differentiable wherever cosine is non-zero. The factor of $\frac{1}{2}$ is the same as the factor seen in the derivatives of the sine and cosine functions. It is nevertheless a curious factor not seen in the Euclidean derivative of the tangent function. 

The derivative of the taxicab cotangent function is found in a similar manner to tangent. The proof will therefore be omitted, but the result is $-\frac{1}{2}\csc^2 \theta$. This again is very reminiscent of the derivative of the Euclidean cotangent function with the added factor of $\frac{1}{2}$.

As we saw earlier, the taxicab secant function is expected to not be differentiable where cosine has an extremum since cosine is not differentiable at these points. This will prevent us from finding a single closed-form expression for the derivative of secant. But, as with tangent we will take the proof in two cases over the intervals $[(k-1)\frac{\pi_t}{2},k\frac{\pi_t}{2})$ with $k$ even and odd. For $k$ even, we have 
\begin{eqnarray*}
\frac{d}{d\theta}(\sec \theta) & = & \frac{d}{d\theta}\left[ \frac{1}{i^{k-2}}\left( k-1-\frac{1}{2}\theta \right)^{-1} \right] \\
        & = & -\frac{1}{i^{k-2}}\left( k-1-\frac{1}{2}\theta\right)^{-2}\left( -\frac{1}{2}\right) \\
        & = & \frac{\frac{1}{2}}{i^{k-2}(k-1-\frac{1}{2}\theta)^2} \\
        & = & \frac{1}{2}\sec \theta \left( \frac{1}{k-1-\frac{1}{2}\theta } \right) \\
        & = & \frac{1}{2}\sec \theta \left( \frac{k-\frac{1}{2}\theta-(k-1-\frac{1}{2}\theta)}{k-1-\frac{1}{2}\theta } \right) \\
        & = & \frac{1}{2}\sec \theta (\tan \theta - 1)
\end{eqnarray*}

For $k$ odd, we have
\begin{eqnarray*}
\frac{d}{d\theta}(\sec \theta) & = & \frac{d}{d\theta}\left[ \frac{1}{i^{k-1}}\left( k-\frac{1}{2}\theta \right)^{-1} \right] \\
        & = & -\frac{1}{i^{k-1}}\left( k-\frac{1}{2}\theta\right)^{-2}\left( -\frac{1}{2}\right) \\
        & = & \frac{\frac{1}{2}}{i^{k-1}(k-\frac{1}{2}\theta)^2} \\
        & = & \frac{1}{2}\sec \theta \left( \frac{1}{k-\frac{1}{2}\theta } \right) \\
        & = & \frac{1}{2}\sec \theta \left( \frac{1-k+\frac{1}{2}\theta+k-\frac{1}{2}\theta)}{k-\frac{1}{2}\theta } \right) \\
        & = & \frac{1}{2}\sec \theta (\tan \theta + 1)
\end{eqnarray*}

This is a very interesting result. Shades of the Euclidean derivative are present, but a whole extra secant term is added or subtracted in the taxicab case. As mentioned above, there is no closed-form expression for the derivative of taxicab secant, and the function has a non-differentiable corner point at $k\pi_t-\frac{\pi_t}{2}$ t-radians.

While this presentation of the taxicab secant derivative contrasts with Euclidean geometry nicely, if we apply the Pythagorean-style trig identities for taxicab geometry the secant derivative becomes even more interesting. When $k$ is even, $\theta$ lives in quadrants II or IV where tangent is negative. In quadrant II, $-\tan \theta +1=-\sec \theta$ and in quadrant IV, $-\tan \theta +1=\sec \theta$. Using these substitutions in the first result above we have the derivative of secant as $\pm \frac{1}{2} \sec^2 \theta$. A similar analysis applies to the second result above. Therefore, on intervals $[k\pi_t,(k+1)\pi_t)$ (excluding the asymptote) with $k$ even and odd,
\[
   \frac{d}{d\theta}(\sec \theta) = \left\{
     \begin{array}{ll}
       \frac{1}{2}\sec^2 \theta, & \text{if }k\text{ is even} \\
       -\frac{1}{2}\sec^2 \theta, & \text{if }k\text{ is odd}
     \end{array}
   \right.
\]
So, unlike the Euclidean secant function, the taxicab secant function is (piecewise) proportional to the square of itself. And, since $\frac{d}{d\theta}\tan \theta = \frac{1}{2}\sec^2 \theta$ the derivative of the secant function is equal to the derivative of the tangent function half of the time!

The derivative of cosecant is found in a similar manner to secant. The proof will therefore be omitted, but the result is shown below. Over the interval $[(k-1)\frac{\pi_t}{2},k\frac{\pi_t}{2})$ with $k$ even and odd, the derivative of cosecant is
\[
   \frac{d}{d\theta}(\csc \theta) = \left\{
     \begin{array}{ll}
       -\frac{1}{2}\csc \theta(\cot \theta - 1), & \text{if }k\text{ is even} \\
       -\frac{1}{2}\csc \theta(\cot \theta + 1), & \text{if }k\text{ is odd}
     \end{array}
   \right.
\]
And similarly, over the interval $[k\pi_t-\frac{\pi_t}{2},k\pi_t+\frac{\pi_t}{2})$ (excluding the asymptote) with $k$ even and odd, the derivative of cosecant can also be expressed as
\[
   \frac{d}{d\theta}(\csc \theta) = \left\{
     \begin{array}{ll}
       -\frac{1}{2}\csc^2 \theta, & \text{if }k\text{ is even} \\
       \frac{1}{2}\csc^2 \theta, & \text{if }k\text{ is odd}
     \end{array}
   \right.
\]

\smallskip

\section{An Alternate Approach}

The derivatives of tangent, cotangent, secant, and cosecant were found above by differentiating combinations of the underlying piecewise functions for sine and cosine. These derivatives can also be found using their basic definitions (such as the ratio of sine and cosine for tangent) and the derivatives of sine and cosine. This approach naturally leads to the secondary results presented above for secant and cosecant. It is also a cleaner approach, but there is benefit in considering the initial results for secant and cosecant found with the presented approach to contrast with the Euclidean derivatives.

\section{Conclusion}

The loss of the derivative relationship between sine and cosine in taxicab geometry has a significant effect on the derivatives of some of the other taxicab trigonometric functions. While the derivatives of the taxicab tangent and cotangent functions retain a very Euclidean flavor (with the seasoning of the factor of $\frac{1}{2}$ from the derivative of sine and cosine), the effect on the secant and cosecant derivatives is to make them (piecewise) proportional to the square of themselves. The derivatives of secant and cosecant are also equal to the derivatives of tangent and cotangent, respectively, half of the time. This is an interesting result, but one that should be expected given the piecewise constant nature of the derivatives of sine and cosine.

\end{document}